\documentclass[
final
]{dmtcs-episciences}

\usepackage[utf8]{inputenc}
\usepackage{subcaption}
\usepackage{amsfonts}
\usepackage{amssymb}
\usepackage{amsmath}
\usepackage{amsthm}
\usepackage{blkarray}
\usepackage{cleveref}
\usepackage{stmaryrd}
\usepackage{tikz}
\usetikzlibrary{positioning}
\usepackage{pgfplots}

\usepackage[round]{natbib}

\newtheorem{theorem}{Theorem}[section]
\newtheorem{lemma}{Lemma}[section]

\newtheorem{openq}{Open Question}[section]
\newtheorem{conjecture}{Conjecture}[section]

\theoremstyle{definition}
\newtheorem{definition}{Definition}[section]

\newtheorem{remark}{Remark}[section]

\author[Thomas C. Hull et al.]{Thomas C. Hull\affiliationmark{1}\thanks{Supported in part by NSF grants DMS-2428771 and DMS-2347000.}
  \and Adham Ibrahim\affiliationmark{2}
  \and Jacob Paltrowitz\affiliationmark{3}\\
  \and Natalya Ter-Saakov\affiliationmark{4}
  \and Grace Wang\affiliationmark{5}}
\title[The stamp folding problem from a mountain-valley perspective]{The stamp folding problem from a mountain-valley
  perspective\thanks{Supported by NSF grant DMS-2149647. The authors thank Mathematical Staircase, Inc. for funding and supporting the MathILy-EST 2024 REU.}}
\affiliation{
  Franklin and Marshall College, Lancaster, USA\\
  Princeton University, Princeton, USA\\
  Harvard University, Cambridge, USA\\
  Rutgers University, New Brunswick, USA\\
  Carnegie Mellon University, Pittsburgh, USA}
\begin{document}
% This is only used if you are compiling for a volume before vol 25
% \publicationdetails{VOL}{2015}{ISS}{NUM}{SUBM}
% This is the new form of collecting the data, starting with vol 25
\publicationdata{vol. 27:3}{2025}{16}{10.46298/dmtcs.15454}{2025-04-01; 2025-04-01; 2025-09-30}{2025-10-03}

%{1998-10-14; 1998-10-14; 2002-07-19; 2014-02-05; 2015-09-09; 2022-12-25}
%{2022-12-3}
%{2022-12-3; None}
%{2023-1-1}
\maketitle
\begin{abstract}
    A strip of square stamps can be folded in many ways such that all of the stamps are stacked in a single pile in the folded state. The stamp folding problem asks for the number of such foldings and has previously been studied extensively. We consider this problem with the additional restriction of fixing the mountain-valley assignment of each crease in the stamp pattern. We provide a closed form for counting the number of legal foldings on specific patterns of mountain-valley assignments, including a surprising appearance of the Catalan numbers. We describe results on  upper and lower bounds for the number of ways to fold a given mountain-valley assignment on the strip of stamps, provide experimental evidence suggesting more general results, and include conjectures and open problems.
\end{abstract}

\section{Introduction}
The stamp folding problem, first formulated by \cite{lucas}, asks for the number of ways to fold a labeled $1 \times n$ strip of stamps into a single flat stack. The total number of ways to fold the $1\times n$ strip, $t(n)$, has been calculated for $n$ up to 45 \cite[A000136]{oeis}. The problem of finding a closed form for the stamp folding numbers is still open, but it is known that $3.065^n<t(n)<4^n$ from \cite{Uehara2016StampFW}.

On the $1 \times n$ strip, we refer to the sections of paper between two creases as faces, and consider these to be labeled in order from left to right, starting from $0$ to $n-1$. We refer to the $n-1$ creases on the strip as tuples defined by the faces they border. So $(i, i+1)$ will refer to the crease bordering faces $i$ and $i+1$. A \textit{mountain-valley assignment} is a function $\mu$ mapping the set of creases to $\{M, V\}$, and refers to the orientation of the creases, so that if we unfold the strip and choose one side to be facing ``up", then mountain (M) creases will look like $\wedge$ and valley creases will look like $\vee$.  In this paper, we consider the question: can we calculate the number of ways to fold the strip of stamps consistent with a particular mountain-valley assignment? Here we consider ways to fold the strip to be different if the faces have different layer orderings (which we make more formal in Section~\ref{sec2}).

This approach to the stamp folding problem is in need of exploration. It is known from  \cite{uehara} and  \cite{UME13} that the \textit{pleat} assignment, the assignment that alternates mountains and valleys, is the only assignment that folds in just one way, and that if every crease is assigned the same direction, then the number of ways to fold is $n-1$. It is noted by   \cite{Uehara2016StampFW} that because there are $2^{n-1}$ MV assignments, the average count must be exponential in $n$, and he provides an explicit example in the conclusion that folds in $\Theta(2^{n/2})$ ways. Mountain-valley assignments have also been considered in connection with the problem of minimizing the crease width, or the number of stamps between each pair of adjacent stamps, from \cite{UME13}. But a general treatment of $1\times n$ stamp-folding problem through a fixed MV assignment lens has not been done.

After formalizing our notation and providing background in Section~\ref{sec2}, we calculate in Section \ref{mv_assign} the number of ways to fold two specific MV assignment patterns, including a new construction of an MV assignment that folds in a number of ways exponential in the length of the string. We conclude with a description of results on bounds, experimental results, conjectures, and open questions.

\section{Background}\label{sec2}

\begin{figure}
    \centering
    \begin{subfigure}{.45\textwidth}
        \centering
        \includegraphics[width=0.8\linewidth]{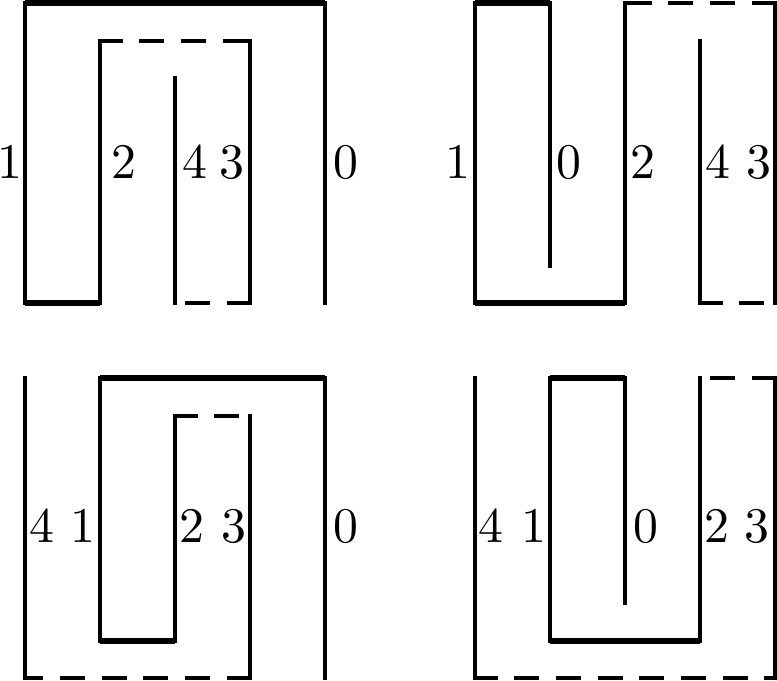}
        \caption{The four possible ways to fold $M^2 V^2$. Face 0 faces to the right. Dashed lines are valley creases.}
        \label{fig:m2v2a}
    \end{subfigure}\hfill
    \begin{subfigure}{.45\textwidth}
        \centering
        \includegraphics[width=\linewidth]{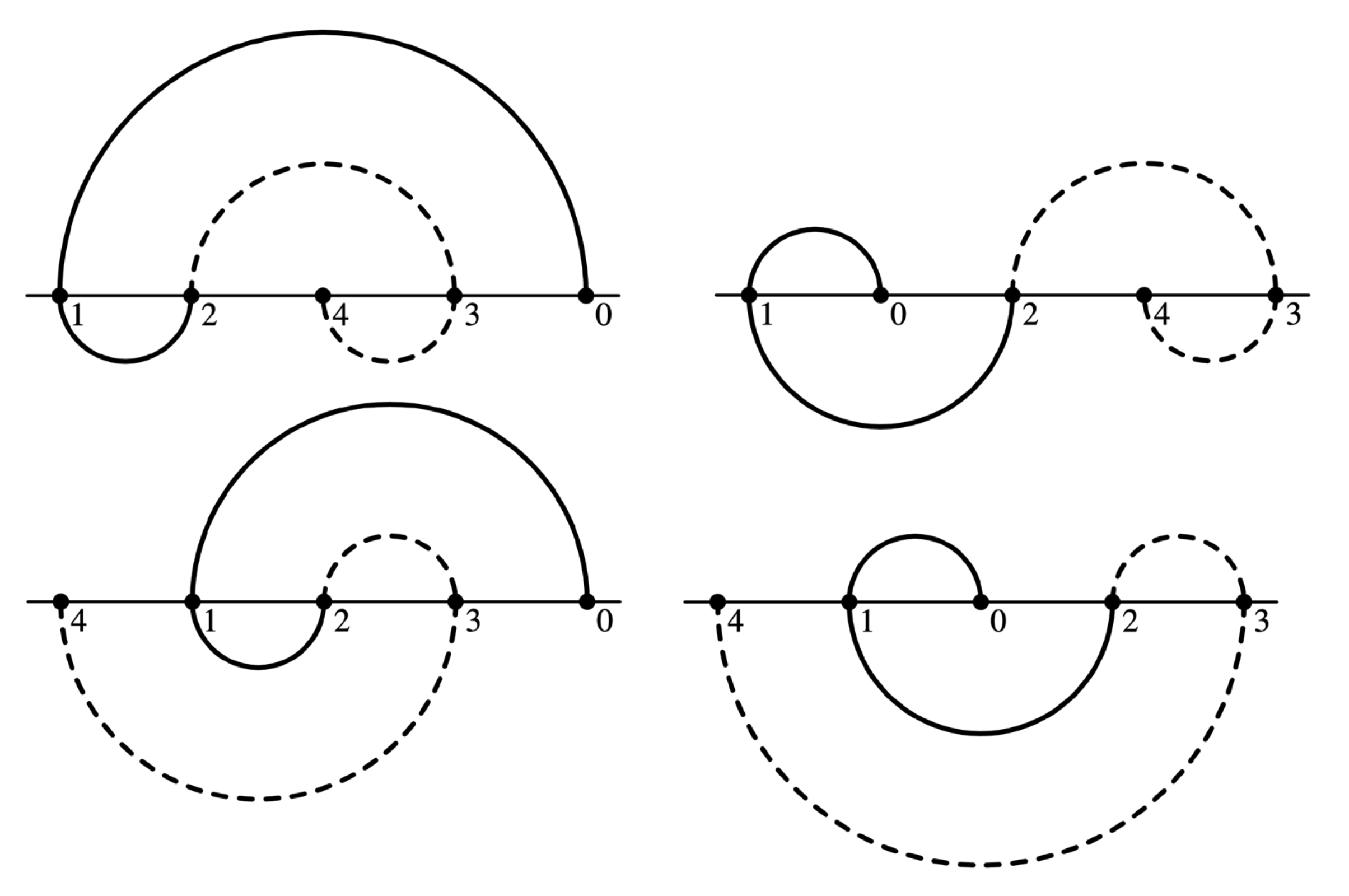}
        \caption{The corresponding meanders. Solid arcs are mountain creases and dashed are valley.}
        \label{fig:m2v2b}
    \end{subfigure}
    \caption{$c(M^2 V^2) = 4$}\label{fig:m2v2}
\end{figure}

We label the faces of our $1\times n$ strip of stamps with numbers $0$ through $n-1$ and consider the labels to all be placed on one side of the paper, so that we may differentiate between the two sides of the paper when folded. We can imagine folding the strip into a flat pile with the faces parallel to one another, as in Figure~\ref{fig:m2v2a} where the faces are drawn vertically. For each crease $(i, i+1)$ in some $1\times n$ strip of stamps, a mountain-valley assignment $\mu$ induces a relation on the adjacent faces $i, i+1$. We say a face is \textit{facing left (right)} if the labeled side of the face is facing that direction in our drawing. Fixing the labeled side face $0$ to face right, we define the \textit{layer ordering} of a folded state to be the permutation of faces from left to right. Fixing face $0$ in this way will then force the other faces to alternate facing left or right, regardless of the MV assignment. That is, all faces labeled an even value will face right as well, and the odd faces will face left. Further, in the folded state, each crease will fold onto either the top or bottom of the vertically drawn cross section in an alternating fashion that is also independent of the MV assignment or the layer ordering. In particular, the crease between $i$ and $i+1$ for even $i$ will fold onto the top side of face $0$, while for odd $i$ it will fold onto the bottom side. See Figure~\ref{fig:m2v2}(a) for examples of a folded $1\times 5$ strip. Using these observations, we can characterize the two conditions that a layer ordering must satisfy.

\begin{definition}\label{valid}
    A permutation $\sigma$ of $\{0, \dots, n-1\}$ is a \textit{valid layer ordering} for a $1 \times n$ with MV assignment $\mu$ if the following hold:
    \begin{enumerate}
        \item (MV relations) Let $i$ and $i+1$ be any two adjacent faces. To respect the MV assignment, if $\mu((i,i+1))=M$ and $i$ is facing right, or  $\mu((i,i+1))=V$ and $i$ is facing left, then $\sigma(i) > \sigma(i+1)$. Otherwise, $\sigma(i) < \sigma(i+1)$.
        %Let $i$ and $i+1$ be any two adjacent faces. To respect the MV assignment, if the crease between them is an $M$ and $i$ is facing right, or the crease is $V$ and $i$ is facing left, then $\sigma(i) > \sigma(i+1)$. Otherwise, $\sigma(i) < \sigma(i+1)$.
        \item (No crossings,  \cite{Koehler}) Let $(i, i+1)$ and $(j, j+1)$ be two creases that fold onto the same side. Then in order for the paper to not self-intersect, it cannot be that $\sigma(i) < \sigma(j) < \sigma(i+1) < \sigma(j+1)$, or any circular permutation of this inequality.
    \end{enumerate}
\end{definition}

For short, we write a MV assignment $\mu$ with the string $\mu((0, 1)) \mu((1, 2)) \dots \mu((n-2, n-1))$. For example, we write the assignment that alternates $M$ and $V$ on the $1 \times 5$ strip as $MVMV$. To shorten strings with many repeated consecutive creases assigned the same value, we denote by $S^k$ the repetition $k$ times of the string $S$. So $(M^2 V^2)^3$ is the same as the string $MMVVMMVVMMVV$. We will refer to a maximal consecutive substring with the same assigned value as a \textit{block}. With this, we describe the counting problem that is the object of study in this paper.

\begin{definition}
    For a given assignment $\mu$, $c(\mu)$ denotes the number of valid layer orderings of the faces with respect to $\mu$. We interchangeably refer to this quantity as $c(S)$, where $S$ is the string associated with $\mu$.
\end{definition}

As an example, Figure \ref{fig:m2v2}(a) shows the four ways that the assignment $M^2 V^2$ on the $1 \times 5$ strip folds, that is, $c(M^2 V^2) = 4$. A related concept that has been previously explored are \textit{meanders}: given a directed non-self intersecting plane curve that intersects a horizontal line $n$ times, what are the permutations that arise from looking at the order of the intersections? For a survey of progress on this problem, see \cite{legendre}. The $1\times n$ stamp folding problem can be seen as a special case where the curve is of finite length: We imagine the stamps being folded vertically as before with the horizontal line drawn through all the faces, and we read off the permutation from left to right. To visualize these, we use solid arcs to represent mountain folds throughout this paper, and dashed arcs to represent valley folds. Figure \ref{fig:m2v2}(b) shows the four meanders representing the ways that $M^2 V^2$ folds. For the remainder of the paper, we will use meanders to represent a folded collection of stamps. 

As further notation for these depictions, we also use $(i, i+1)$ to refer to the arc in a meander drawing between faces $i$ and $i+1$, which corresponds to the crease between these faces. It will also be useful to refer to intervals along the horizontal line from the meander interpretation (which we will simply call \textit{the line}), so we let $I(i, j)$ be the interval bounded by $i$ and $j$ on the line (so $I(i, j)$ and $I(j, i)$ refer to the same interval). In this way, we are abstracting away from physical paper, which would require that arcs between consecutive faces be of the same length, and focusing only on the relative order of faces.

Our proofs will rely on the following core idea.

\begin{lemma} \label{prop:induction}
    Let $\sigma$ be a valid layer ordering with respect to a MV assignment $\mu$ on a $1 \times n$ strip of stamps. Define $\mu'$ on a $1 \times (n-1)$ strip to be the assignment identical to $\mu$ for all $(i, i+1)$ with $i < n-1$. Then, the permutation $\sigma'$ resulting from removing face $n-1$ from $\sigma$ gives a valid layer ordering with respect to $\mu'$.
\end{lemma}

The proof is a direct application of the definition of a valid layer ordering. By induction, this means that we can recursively build up the valid layer orderings for any MV assignment one face at a time: After having placed $0$ through $i-1$ into a permutation that is a valid layer ordering for the first $i$ faces, we look for possible positions to place the next face that will maintain the MV relations and not induce crossing arcs. With the meander visualization, this corresponds to finding an interval $I(j, k)$ for $j, k < i$ that is on the appropriate side of $i-1$ (for example, if $\sigma(i)$ must be greater than $\sigma(i-1)$ we look for $j, k$ with $\sigma(j), \sigma(k) > \sigma(i-1)$), where the arcs previously drawn will not intersect the arc $(i-1, i)$.

We can begin with an observation on the counts for related MV assignments. 

\begin{remark}\label{symm}
    The number of valid layer orderings of an MV assignment is invariant under reversing the string corresponding to the MV assignment, as well as flipping all of the mountains to valleys. 
\end{remark}

This equivalence is due to the fact that these MV assignments correspond to the same folded piece of paper, just rotated or reflected. As an example of using the perspective of meanders, we demonstrate a proof of a known simple counting result, first stated in \cite{uehara} without proof.

\begin{theorem}\label{m^n}
    For every $n$, on the $1 \times (n+1)$ strip we have $c(M^n) = c(V^n) = n$.
\end{theorem}

\begin{proof}
We will prove this for the MV assignment $M^n$, which is equivalent to $V^n$ by Remark \ref{symm}. For this MV assignment, the folding consists of two spirals, one on each end of the strip, as seen in Figure~\ref{fig:prelim1}. An intuitive approach is to simply note that the $n+1$ faces can be partitioned into two sets, $n+1=j+k$, in $n$ ways, but this does not address whether such a partition could be realized in more than one way. To handle this and prepare the reader for forthcoming proofs, we formalize the argument by considering how to add an additional face to an existing folding. Suppose we have an existing valid layer ordering for $n$ stamps using the MV assignment $M^{n-1}$. In order to extend this to a valid layer ordering on $n+1$ stamps with the MV assignment $M^{n}$, we must consider two cases:

\begin{enumerate}[{Case }1:]
\item The $n$-th face is not on the ``outside" of the folded strip, i.e., $\sigma(n)\notin \{0, n\}$ in the existing layer ordering for $M^{n-1}$. An example of this case is seen in Figure~\ref{fig:prelim1}(a). Then, we can see inductively that face $n$ must lie in $I(n-1, n-2)$, forcing face $n+1$ to continue the spiral and lie in $I(n-1, n)$.

\begin{figure}
\centerline{\includegraphics[width=\linewidth]{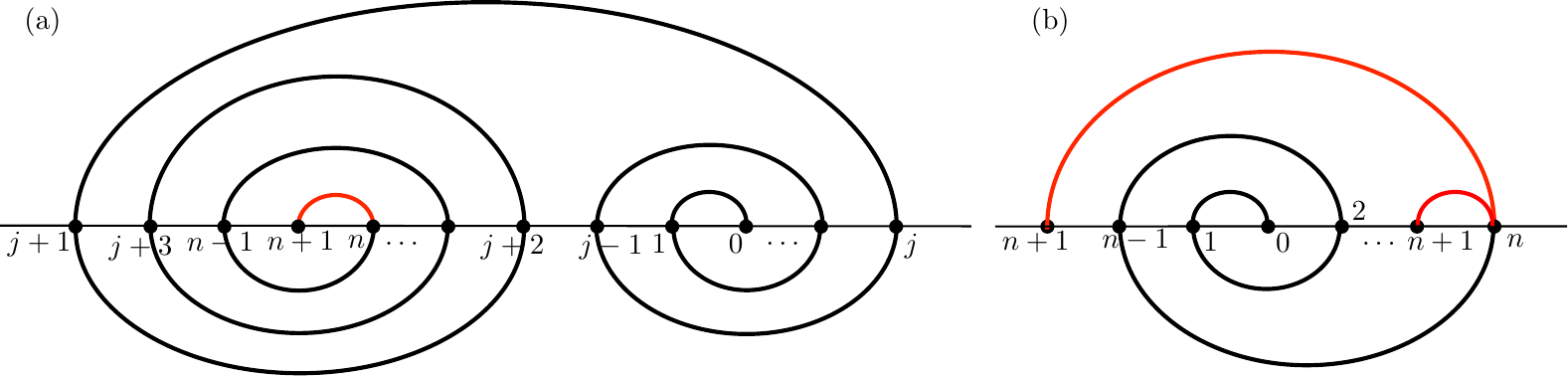}}
        \caption{(a) Addition of a single face in Case 1 and (b) Case 2 of Theorem~\ref{m^n}.}
        \label{fig:prelim1}
\end{figure}

\item The $n$-th face is on the outside of the folded strip in the existing layer ordering for $M^{n-1}$. As seen in Figure \ref{fig:prelim1}(b), this case allows for two places to put stamp $(n+1)$, namely $I(n-2,n)$ or $I(\pm\infty,n-1)$, with exactly one of the ways keeping the last face on the outside.
\end{enumerate}

Thus we have that $c(M^n) = c(M^{n-1})+1$.  Solving this recurrence with the initial condition $c(M^1) = 1$ gives us $c(M^n)=n$.
 \end{proof}
    
\section{Specific Mountain-Valley Assignments}\label{mv_assign}

\subsection{Two-Block  Assignments}\label{ab_pattern}
Let's count the number of valid layer orderings in MV assignments with only two blocks. Removing symmetries, it suffices to consider MV assignments of the form $M^a V^b$ with $a\geq b$ by Remark~\ref{symm}. 

\begin{theorem}\label{thm:MaVb}
    The MV assignment $M^aV^b$ on the $1 \times (a+b+1)$ strip, where $a > b + 1, b> 4$ folds in $ab+b^2-b+1$ ways. 
\end{theorem}

\begin{proof}
Beginning with one of the $a$ valid stamp foldings for $M^a$, we will calculate the number of ways to fold the second block, $V^b$. Let the extreme left and right
faces of this folding of $M^a$ be $j$ and $j+1$, for some $j \in {0, \dots, a-1}$. 

If $j+1=a$ then the folding of $M^a$ is a single spiral starting with face 0 in the center, as seen in Figure~\ref{fig:MaVb}(a). In this case face $a+1$ must be folded in $I(a,\pm\infty)$ and we have choices for the remaining $b-1$ faces. Face $a+2$ can go into $I(a,a+1)$, in which case it must deterministically spiral inwards from there. In fact, this option of deterministically spiraling inwards is always an option at every decision point in folding $V^b$. The other options are to spiral outwards, shown in red in Figure~\ref{fig:MaVb}(a) (i.e., fold into $I(a-1,\pm\infty)$, then into $I(a+1,\pm\infty)$, and so on), or to spiral into the $M^a$ spiral by folding into $I(a-1,\pm\infty)$, then into $I(a-2,a)$ and following the blue path in Figure~\ref{fig:MaVb}(a). Including the options of spiraling deterministically at every juncture, there are $b$ options in the outward-spiraling red paths and $b-2$ for the inward-spiraling blue paths. Thus the $j+1=a$ case has $2b-2$ ways to fold.

\begin{figure}
    \centering
    \includegraphics[width = \linewidth]{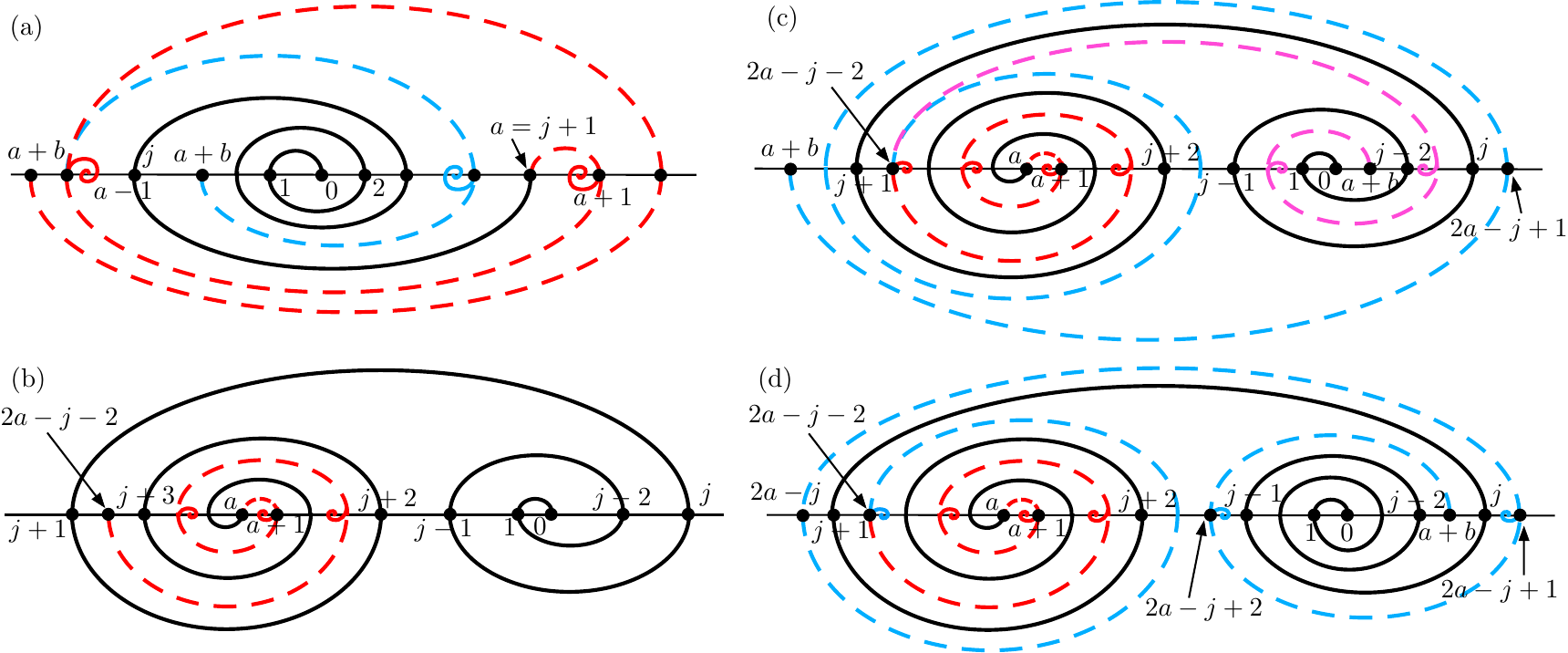}
    \caption{The various cases of Theorem~\ref{thm:MaVb}.}
    \label{fig:MaVb}
\end{figure}

If $j+1<a$ then, given that $a>b+1>4$, the general situation will look something like that shown in Figure~\ref{fig:MaVb}(b), where the $V^b$ folding spirals away from stamp $a$ (shown in red) and assuming it does not deterministically spiral inwards at some point, will reach stamp $2a-j-2$. From there, our folding of $V^b$ has three options: 

\begin{enumerate}[(1)]
\item spiral deterministically at that point, (the small red spiral at stamp $2a-j-2$ in Figure~\ref{fig:MaVb}(c)), 

\item intermix the valley creases with the mountain creases in the other spiral by entering $I(j-2,j)$ (or $I(j-1,j)$ if $j=1$), shown in pink in Figure~\ref{fig:MaVb}(c), or 

\item not intermix the valley creases by entering $I(j-1,j+2)$ (or $I(j,j+2)$ if $j=0$) and spiral or circle outward, shown as blue paths in Figure~\ref{fig:MaVb}(c) and (d).
\end{enumerate}

The only options that do not depend on $j$ are  (1) and (3) when the folding of $V^b$ spirals away outside the folding of $M^a$. We consider these together, with the red path (and its deterministic spiral options at each juncture) in Figure~\ref{fig:MaVb}(b) and the blue path in Figure~\ref{fig:MaVb}(c). This has $b$ ways to fold, as there are $b-1$ places to spiral inwards (including the option (1) case) plus the one case where the whole path spirals outward. Any of these $b$ ways to fold $V^b$ can happen with any of the $a$ ways to fold $M^a$ except for the $j+1=a$ case described earlier. Thus there are $(a-1)b=ab-b$ ways to fold this case.

Option (2) assumes we have not previously spiraled deterministically at a previous juncture and only happens if $2a-j-1\leq a+b$ (otherwise the red path in Figure~\ref{fig:MaVb}(b) never leaves the spiral centered at stamp $a$). That is, $a-b-1\leq j$. If $j=a-b-1$ then we have only one way to fold the one remaining stamp of $V^b$. With each increase in $j$ we have one more stamp of $V^b$ to fold along the pink path of Figure~\ref{fig:MaVb}(c), with the extreme case being $j=a-2$ and all $b$ stamps of $V^b$ folding along the pink path, which has $b$ ways to fold. Thus in this case, ranging over all possible values of $j$, we have 
$1+2+\dots+b = b(b+1)/2$ ways to fold in this case.

The remaining case is option (3) where our folding enters $I(j-1,j+2)$, loops around the $M^a$ folding, and then stamp $2a-j+2$ re-enters at $I(j-1,2a-j-1)$, as shown in the blue path in Figure~\ref{fig:MaVb}(d). For this to occur, we need $2a-j+2\leq a+b$, or $a-b+2 \leq j$. Note that $j$ is at most $a-2$, so for this case to exist, $b$ must be at least 4. This is essentially the same as the previous case, except we make three folds before entering $I(j-1,2a-j-1)$. Thus there are $1+2+\dots+(b-3)=(b-2)(b-3)/2$ ways to fold this case.

Putting all these cases together, for $a - b > 1$ and $b\geq 4$,  there are 
$$    (2b-2)+(ab - b) + \frac{b(b+1)}{2} + \frac{(b-2)(b-3)}{2} = ab +b^2 - b +1$$
ways to fold the MV assignment $M^aV^b$. 
\end{proof}

\subsubsection{Edge Cases}

There are also the edge cases that need to be considered. First, if $b< 4$, the last case in the above proof cannot occur, and thus the $\frac{(b-2)(b-3)}{2}$ must be removed from the sum. Note that this only affects the sum when $b < 2$, as $\frac{(b-2)(b-3)}{2}$ equals zero when $b$ is either $2$ or $3$. 

Also,  if $a-b$ is zero or one, the folding path can reach the center of the second spiral, which limits the number of elements of the option (2). So, in option (2) where $a-b = 1$, we only have $1+2+\dots + b-1 = b(b-1)/2$ paths, and $M^a V^{a+1}$ folds in $2a^2-a+1$ ways.

In the case where $a = b$, we have two fewer steps until reaching the center of the spiral, and thus option (2) has $1+2+\dots+b-2$ elements, making $M^a V^a$ fold in $2a^2-3a+2$ ways. 

\begin{remark}\label{rem:poly}
    Theorem~\ref{thm:MaVb} indicates that if $n=a+b$ is the number of creases in our strip of stamps, then $c(M^a V^b)$ is polynomial in $n$ and is $O(n^2)$ in particular. Intuitively, the reason this two-block case is polynomial is due to the fact that if an additional stamp is added, either in the M or the V block, the number of intervals into which it can be inserted is bounded by a constant; it will have either one, two, or three choices, as seen in the proof. In the next section, we will see MV assignments whose number of ways to fold is exponential in $n$. 
\end{remark}

\subsection{The 2-alternating Assignments}\label{catalan_pattern}

We have seen a class of MV assignments that fold in a number of ways polynomial in the length of the string. In this section we consider the class of assignments with blocks of size two. The first few examples are $MM$, $MMVV$, $MMVVMM$, and so on. We refer to these assignments as \textit{2-alternating}. Formally, the 2-alternating assignment with $m$ blocks is $(M^2 V^2)^{m/2}$ for even $m$ and $(M^2 V^2)^{(m-1)/2} M^2$ for odd $m$. We now aim to count precisely the number of ways this pattern folds.

\begin{theorem}\label{thm:m2v2}
    The 2-alternating assignment on the $1 \times (2m+1)$ strip with $m$ blocks folds in $\displaystyle C_{\lfloor \frac{m}{2} \rfloor + 1} C_{\lceil \frac{m}{2} \rceil + 1}$ ways, where $C_i$ is the $i$-th Catalan number. 
\end{theorem}

The proof has two main parts. First, we show that counting the valid layer orderings for these MV assignments is equivalent to counting restricted walks in $\mathbb{Z}^2$ in Lemma \ref{lem:walks}. We then show that the number of such walks is the product of Catalan numbers as specified in the Theorem.

\begin{figure}
\centerline{\includegraphics[width=\linewidth]{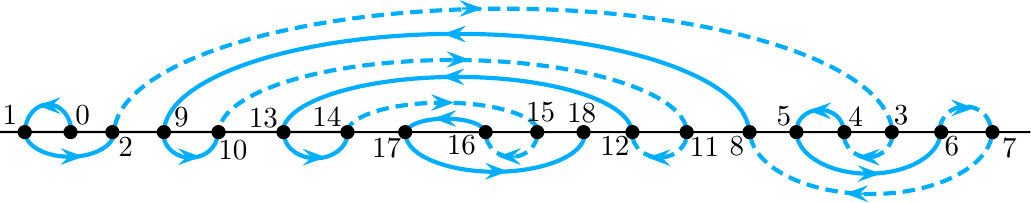}}
\caption{An example of a 2-alternating assignment folding of $(M^2V^2)^4 M^2$.}\label{fig:M2V2ex}
\end{figure}

To build up to this, we make some observations about the inherent structure of 2-alternating assignments; an example is shown in Figure~\ref{fig:M2V2ex}. Observe that when viewing our foldings as meanders, arcs $(i, i+1)$ lie above the horizontal line for even $i$, and below for odd $i$; we refer to these as \textit{top arcs} and \textit{bottom arcs} respectively. Note that our convention is that the first face $0$ faces right and the first $M$ crease $(0,1)$ is a top arc, which means in the meander drawing this arc will be drawn from right-to-left; i.e., $\sigma(1)<\sigma(0)$. Perpetuating this, we have that each initial $M$ in an $M^2$ will be a top arc traveling right-to-left in the meander drawing, while the second $M$ will be a bottom arc going left-to-right. The $V^2$ instances will be similar but with the directions reversed. 

 We call a bottom arc $(i, i+1)$ \textit{uncovered} if $I(i, i+1)$ is not contained in any $I(j, j+1)$ for bottom arc $(j, j+1)$, with $j < i$. Otherwise, we call $(i, i+1)$ \textit{covered}. These special arcs will be key to counting the valid layer orderings.

\begin{lemma}\label{lem:m-first-bot}
    In every possible layer ordering satisfying a 2-alternating assignment, the two following conditions are satisfied:
    \begin{enumerate}
        \item Suppose $(i, i+1)$ is a top arc in a 2-alternating assignment. Then given a layer ordering of faces 0 through $i$, face $i+1$ has precisely one place to be inserted in the permutation.
        \item The uncovered bottom arcs are arranged such that all of the mountain arcs are to the left of all the valley arcs. That is, if $(i, i+1)$ is an uncovered bottom mountain arc and $(j, j+1)$ is an uncovered bottom valley arc, then $\sigma(i) < \sigma(j)$. Further, if the assignment ends with a $V$, the position of the final face must be the endpoint of the leftmost uncovered bottom valley arc. Otherwise, the position of the final face must be the endpoint of the rightmost uncovered bottom mountain arc.
    \end{enumerate}
\end{lemma}

\begin{figure}
    \centering
    \includegraphics[width=0.8\textwidth]{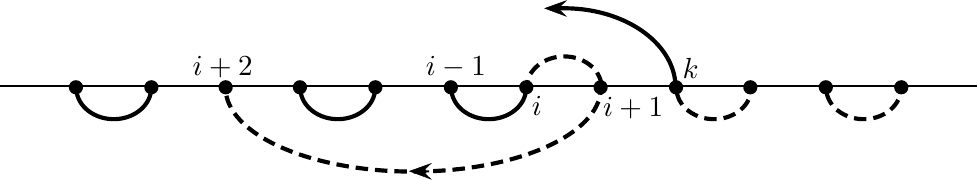}
    \caption{Illustration of the proof of Lemma \ref{lem:m-first-bot}. Not all arcs are shown.}
    \label{fig:main-lemma}
\end{figure}

\begin{proof} 

    We prove this by induction on the number of stamps in a 2-alternating assignment. It is trivially true for the first bottom mountain arc. Now suppose the lemma is true for 2-alternating assignments on $i$ stamps. Given a valid folding of a 2-alternating assignment on $i+1$ stamps, drawn as a meander, we have that our meander picture satisfies the lemma for points $0$ up to $i$  with $(i-1,i)$ being the endpoint of a bottom arc, which we for the moment assume is a mountain arc. See Figure \ref{fig:main-lemma} for an illustration of this.
   
Let $k$ be the left-most face in our drawing that is part of the left-most bottom valley arc. Then by the induction hypothesis we have $\sigma(i)< \sigma(k)$.
    Next, we will show that $\sigma(i+1) < \sigma(k)$ and thus there is only one place for face $i+1$ to be folded: $I(i, k)$ (which will prove part (1) of the lemma by induction). Suppose for the sake of contradiction that $\sigma(k) < \sigma(i+1)$. Since $(k-1, k)$ is a valley by our inductive hypothesis, the arc $(k, k+1)$ must be a top $M$ arc with $\sigma(k+1) < \sigma(k)$. We also know that $\sigma(k+1) < \sigma(i)$. If not, $(k+1,k+2)$ would be a bottom mountain arc contained in $I(i, k)$ violating the assumption of the inductive hypothesis that $(i-1, i)$ is the rightmost mountain arc. But then $\sigma(k+1) < \sigma(i) < \sigma(k) < \sigma(i+1)$, so $(i,i+1)$ and $(k,k+1)$ intersect and the layer ordering is not valid. 

    Then $(i+1, i+2)$ will be a bottom $M$ arc in the left direction that folds face $i+2$ into an interval to the left of $I(j, k)$, which maintains property (2) of the lemma. The case where $i$ was the endpoint of a valley arc is symmetric and the same reasoning applies.
\end{proof}

It is important to note that this structure would not be present if the MV assignment was different. For general assignments, the uncovered arcs could be intermixed between mountain and valley, and it is rare that the structure would require that any face has only one place to fold.

\begin{figure}
    \centering
    \includegraphics[width=0.8\textwidth]{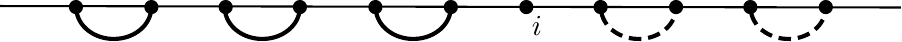}
    \caption{With only uncovered bottom arcs shown, at $i$ there are $4$ potential gaps to the left and $3$ to the right. If $(i, i+1)$ is a mountain arc then $f(i) = 4$ and $b(i) = 3$. Otherwise $f(i) = 3$ and $b(i) = 4$. }
    \label{fig:fb}
\end{figure}

This means that to count the total ways to build the folding, we need only pay attention to the bottom arcs and where they go. The utility of Lemma \ref{lem:m-first-bot} is that we need only maintain the number of uncovered bottom arcs of both kinds. Consider an example illustrated in Figure \ref{fig:fb}: Say we have already folded the faces up to $i$ and there are three uncovered bottom mountain arcs and two valley arcs. Then face $i+1$ with a bottom arc $(i, i+1)$, must fold in the gaps (uncovered intervals) between these arcs or the gap to the far left or right. If we know that $(i, i+1)$ is a valley, there are four options including entering the left-most, unbounded gap before the rightmost mountain.  Otherwise, if $(i, i+1)$ is a mountain, there are three options.

To capture this idea, define $f(i)$ to be the number of potential gaps for $(i, i+1)$ to enter  in the ``forward'' direction; if $(i, i+1)$ is a valley arc then these are the gaps to left of $i$, otherwise they are those to the right. Then define $b(i)$ to be the gaps ``behind'' $i$. Notice that between consecutive bottom arcs, the direction of forwards and backwards switches. We then have the following.

\begin{lemma}\label{lem:recurrence}
    Let $i$ be the starting point of a bottom arc with values $f(i)$ and $b(i)$. Then, $f(i+2) = b(i) + 1$ and $b(i+2) \in \{1, \dots, f(i) \}$, where $b(i+2)$ depends on where face $i+1$ is folded.
\end{lemma}

\begin{proof}
    Suppose face $i+1$ is placed between $j$ and $k$. Then for $i+2$, the gaps in the direction of the arc will be those gaps that were behind $(i, i+1)$, including the newly created gap $I(i+1, k)$. Conversely, because $i+1$ can be folded into any of the $f(i)$ gaps and it will create a new one in front of it, this will give some number of gaps behind $(i+2, i+3)$ that is at minimum 1 and at most $f(i)$.
\end{proof}

We can now completely abstract the process of counting valid layer orderings to the following.

\begin{lemma}\label{lem:walks}
    Let $\mu$ be the 2-alternating assignment on a $1 \times (2m+1)$ strip. Then $c(\mu)$ is given by the number of walks of length $m$ in $\mathbb{Z}^2$ starting at $(2, 1)$ with the allowed steps of the form $(a, b) \to (b+1, i)$ where $i \in \{1, \dots, a\}$.
\end{lemma}

\begin{proof}
    Starting with $i=1$, we let $f(i)=f(1) = 2$ and $b(i) =b(1)= 2$ and imagine a walk in $\mathbb{Z}^2$ as described starting at $(2, 1)$. At each step in our walk we add 2 to the index $i$ and let $(f(i+2), b(i+2))$ be point in $\mathbb{Z}^2$ in our walk. Thus our walk can bring us to any of the positions described in Lemma \ref{lem:recurrence}, and any of these walks in $\mathbb{Z}^2$ can determine an $M^2$ or a $V^2$ part of our 2-alternating stamp folding. We take $m$ steps to reach index $2m+1$, the final index in the permutation which we know to be of length $2m+1$.
\end{proof}

To count the number of ways to perform this walk, we introduce a matrix that will encode the state information at each step.

\begin{definition}\label{d_X}
    Let $i,j\in\mathbb{N}$. We define a matrix $X(m)$ such that $X(m)_{i-1, j}$ is the number of walks of length $m$ in $\mathbb{Z}^2$ starting at $(2,1)$ with steps as defined above that end at $(i, j)$. We omit rows and columns of all zeroes. See below for an example when $m=4$.
\end{definition}
$$
      X(4) = \begin{blockarray}{ccccll}
        & & j \\
        & \textbf{1} & \textbf{2} & \textbf{3} & & \\
        \begin{block}{[cccc]ll}
          & 4 &  4 & 2 && \textbf{2} \\
          & 4 & 4 & 2 && \textbf{3} & $i$ \\
          & 2 & 2 & 1 && \textbf{4} \\
        \end{block}
      \end{blockarray}
$$

One key feature we want to know about $X(m)$ is its dimensions. This will be important later as we construct a matrix representation of taking a single step that will rely on the dimensions of $X(m)$.

\begin{lemma}\label{l_s1}
    Let $(a,b)$ be the position in a walk after $m$ steps. Then $2\leq a\leq \lfloor \frac{m+4}{2}\rfloor$ and $1\leq b\leq \lceil\frac{m+2}{2}\rceil$.
    \end{lemma}
    
    \begin{proof}
    We proceed by induction on $m$. Consider $(a,b)$ a position in the walk after $m+1$ steps. Then by definition, $b\geq 1$. Also by definition, there is an element $(c,d)$ in the walk after $m$ steps such that $a=d+1$ and $b\leq c$. By the inductive hypothesis, $d\geq 1$. So $a\geq 2$. Also by the inductive hypothesis, $c\leq \lfloor \frac{m+4}{2}\rfloor $ and  $d\leq \lceil\frac{m+2}{2}\rceil$. Thus $b\leq \lfloor \frac{m+4}{2}\rfloor =\lceil \frac{(m+1)+2}{2}\rceil$ and  $a\leq \lceil\frac{m+4}{2}\rceil=\lfloor\frac{(m+1)+4}{2}\rfloor$.
    \end{proof}

Now, because of Lemma~\ref{l_s1}, we can see that $X(m)$ must always have dimension $\lfloor \frac{m+2}{2}\rfloor\times \lceil\frac{m+2}{2}\rceil$. We say that accessing any index outside of these bounds will return a value of 0. We can now define the matrix representation of taking a single step in a walk on $\mathbb{Z}^2$ mentioned above.

\begin{definition}\label{d_A}
    Define the matrix $A(m)$ to be the $(\lfloor \frac{m+2}{2}\rfloor+1)\times \lfloor \frac{m+2}{2}\rfloor$ matrix $$\begin{bmatrix}
        1 & 1 & ... & 1 & 1 \\
        1 & 1 & ... & 1 & 1 \\
        0 & 1 & ... & 1 & 1 \\
        & & ... & & \\
        0 & 0 & ... & 1 & 1 \\
        0 & 0 & ... & 0 & 1
    \end{bmatrix}$$
    That is, $A(m)$ is an upper triangular matrix of all ones with a row of ones appended to the top.
\end{definition}

Now we prove that the definition of $A(m)$ given above is actually helpful in taking a step.

\begin{lemma}\label{l_AXT}
    Applying $A(m)$ to $X(m)$ and then transposing the result is the same as taking a step from the positions represented by $X(m)$. In other words,
    $$X(m+1) = (A(m)X(m))^T$$
\end{lemma}
\begin{proof}
    When taking a step at position $(c, d)$, for each positive integer $i$ less than or equal to $c$, we can move to position $(d+1, i)$. Let $(a,b)$ be a position in the walk that can be reached after $m+1$ steps. The number of ways to reach $(a,b)$ in $m+1$ steps is thus equal to the sum of the number of ways to reach each element $(i, a-1)$ in $m$ steps, where $i\geq b$. Since $X(m+1)_{a-1, b}$ is the number of ways to reach position $(a,b)$ in $m+1$ steps, $$X(m+1)_{a-1, b}=\sum_{i=b-1}^{\lfloor\frac{m+2}{2}\rfloor}X(m)_{i, a-1}$$
    This operation is the same as applying $A(m)$ to $X(m)$ and then transposing the result.
\end{proof}

Now that we have a matrix for both taking a single step and representing all valid positions after taking $m$ steps, we can use matrix multiplication to easily compute the entries in $X(m)$ for any $m$. 

\begin{figure}[h]
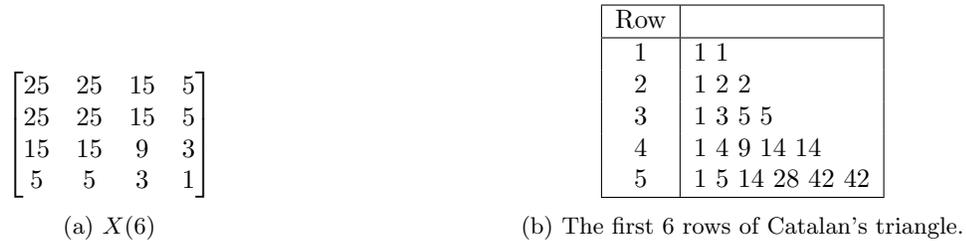

    \centering
    \begin{subfigure}{0.49\textwidth}
        \centering
        $\begin{bmatrix}
        25 & 25 & 15 & 5 \\
        25 & 25 & 15 & 5 \\
        15 & 15 & 9 & 3 \\
        5 & 5 & 3 & 1
        \end{bmatrix}$
        \caption{$X(6)$}
        \label{X(6)}
    \end{subfigure}
    \hfill
    \begin{subfigure}{0.49\textwidth}
        \centering
        \begin{tabular}{ | c | l | }
          \hline
          Row & \\
          \hline
          1 & 1 1 \\
          2 & 1 2 2 \\
          3 & 1 3 5 5 \\
          4 & 1 4 9 14 14 \\
          5 & 1 5 14 28 42 42 \\
          \hline
        \end{tabular}
        \caption{The first 6 rows of Catalan's triangle.}
        \label{cat_tri}
    \end{subfigure}
    \caption{The third row of Catalan's Triangle flipped backwards is equal to the last row and column of $X(6)$.}\label{Xcat}
\end{figure}

\begin{definition}
    \textit{Catalan's Triangle} is a triangle of numbers such that each entry is equal to the sum of the entries of the previous row that are not to the right of the current position. One important property of Catalan's triangle is that the elements in the $n$-th row sum to the $n$-th Catalan number. The first five rows of Catalan's triangle are shown in \Cref{cat_tri}. For more on this triangle, see \cite{doi:10.1080/0025570X.1996.11996408}.  
\end{definition}

It turns out we can connect the entries of our matrices $X(m)$ to this triangle. (See \Cref{Xcat}(a) and \Cref{Xcat}(b) for an example.)

\begin{lemma}\label{l_rcol}
    The last column of $X(m)$ is equal to the $\lfloor\frac{m}{2}\rfloor$-th row of Catalan's triangle and the last row of $X(m)$ is equal to the $\lceil\frac{m}{2}\rceil$-th row of Catalan's triangle. 
\end{lemma}

\begin{proof}
    We proceed by induction on $m$. $X(2)$ is the $2\times 2$ matrix of all ones, whose last row and column equals the first row of Catalan's triangle.
    
    By \Cref{l_AXT}, $X(m+1)=(A(m)X(m))^T$. We will first prove the claim holds for the last column. Let $i\in\mathbb{N}$ such that $1\leq i\leq \lfloor \frac{(m+1)+2}{2}\rfloor$. By matrix multiplication, 
    $$X(m+1)_{i, \lceil\frac{(m+1)+2}{2}\rceil} = \sum_{j=\lceil\frac{m+1}{2}\rceil}^{\lfloor\frac{m+2}{2}\rfloor}X(m)_{j, i}= X(m)_{\lfloor\frac{m+2}{2}\rfloor, i}$$
    So the last column of $X(m+1)$ is equal to the last row of $X(m)$, which is the $\lceil\frac{m}{2}\rceil=\lfloor\frac{m+1}{2}\rfloor$th column of Catalan's triangle by the inductive hypothesis.
    
    Now we prove the claim holds for the last row. Let $i\in\mathbb{N}$ such that $1\leq i\leq \lceil \frac{(m+1)+2}{2}\rceil$. By matrix multiplication, 
    $$X(m+1)_{\lfloor\frac{(m+1)+2}{2}\rfloor, i} = \sum_{j=i-1}^{\lceil\frac{m+2}{2}\rceil}X(m)_{\lfloor\frac{m+2}{2}\rfloor, j}$$
    So the last column of $X(m+1)$ is dependent only on the last row of $X(m)$, which by our inductive hypothesis is equal to the $\lceil\frac{m}{2}\rceil$-th row of Catalan's triangle, and this recurrence is identical to the one that generates the Catalan triangle.
\end{proof}

We now know the last row and column of any $X(m)$, but we still need to know the other entries. Interestingly,  each element of the matrix is equal to the product of the element in its same row and last column and the element in its same column but last row. This will be defined rigorously in the following lemma.

\begin{lemma}\label{l_matin}
    Let $i,j\in\mathbb{N}$ such that $1\leq i\leq \lfloor \frac{m}{2}\rfloor$ and $1\leq j\leq \lceil \frac{m}{2}\rceil$. Then $$X(m)_{i, j} = X(m)_{i, \lceil\frac{m+2}{2}\rceil}\cdot X(m)_{\lfloor \frac{m+2}{2}\rfloor, j}$$
\end{lemma}

\begin{proof}
    We proceed by induction on $m$. By \Cref{l_AXT}, $X(m+1)=(A(m)X(m))^T$. $A(m)$ has dimension $(\lfloor \frac{m+2}{2}\rfloor+1)\times \lfloor \frac{m+2}{2}\rfloor$ by \Cref{d_A}. Let $i,j$ be natural numbers such that $1\leq i\leq \lfloor \frac{m+1}{2}\rfloor$ and $1\leq j\leq \lceil \frac{m+1}{2}\rceil$. Using matrix multiplication, we see that
    \begin{align}
        X(m+1)_{i, j} &= \sum^{\lfloor \frac{m+1}{2}\rfloor}_{k=i-1}\sum^{\lceil \frac{m+1}{2}\rceil}_{\ell=j-1}X(m)_{k, \ell}\notag\\
        &= \sum^{\lfloor \frac{m+1}{2}\rfloor}_{k=i-1}X(m)_{k, \lceil \frac{m+1}{2}\rceil}\sum^{\lceil \frac{m+1}{2}\rceil}_{\ell=j-1}X(m)_{\lfloor \frac{m+1}{2}\rfloor, \ell}\tag{Inductive Hypothesis}
    \end{align}
    By the Catalan's triangle recurrence and \Cref{l_rcol}, this is equal to $X(m+1)_{i, \lceil\frac{(m+1)+2}{2}\rceil}\cdot X(m+1)_{\lfloor \frac{(m+1)+2}{2}\rfloor, j}$.
\end{proof} 

\begin{proof}[of \Cref{thm:m2v2}]
    The number of walks of length $m$ in $\mathbb{Z}^2$ satisfying the above conditions is equal to the sum of the elements in $X(m)$ since $X(m)$ contains the multiplicity of each reachable position after $m$ steps. By \Cref{l_matin}, the sum of the elements in $X(m)$ is equal to the sum of the elements in its last column multiplied with the sum of the elements in its last row. By \Cref{l_rcol}, these elements are exactly the elements in the $\lfloor\frac{m}{2}\rfloor$-th and $\lceil\frac{m}{2}\rceil$-th rows of Catalan's triangle, respectively. So the sum of the elements in $X(m)$ is equal to the sum of the elements in the $\lfloor\frac{m}{2}\rfloor$-th row of Catalan's triangle multiplied by the elements in the $\lceil\frac{m}{2}\rceil$-th row, which are known to be the $\lfloor\frac{m}{2}\rfloor+1$-st and $\lceil\frac{m}{2}\rceil+1$-st Catalan numbers, respectively.
\end{proof}

\begin{remark}
    Theorem~\ref{thm:m2v2} proves that the number of ways to fold 2-alternating assignments $(M^2V^2)^n$ and $(M^2V^2)^n M^2$ increases exponentially with $n$. Intuitively, this is because as we add more stamps to our 2-alternating assignment there can be an increasing number of intervals in which to place the next stamp, as captured by the $f(i)$ and $b(i)$ functions of Lemma~\ref{lem:recurrence}, whereas the polynomial growth of Theorem~\ref{thm:MaVb} never had more than three available intervals for the next stamp. 
\end{remark}

\section{Conclusion}

We have given a closed form for the number of ways to fold two specific MV assignment families for a $1\times n$ strip of stamps, one polynomial and one exponential in the number of stamps. In a previous, expanded version of this paper (\cite{stamp1}), we prove non-tight bounds on general fold enumerations for a given MV assignment. Specifically, if we let $\llbracket a_1, a_2, \dots, a_m\rrbracket$ refer to the MV assignment $M^{a_1}V^{a_2}\dots (M/V)^{a_m}$, then we can prove the following:

\begin{theorem}[\cite{stamp1}]
    For an arbitrary mountain-valley assignment $\llbracket a_1, a_2, \dots, a_m\rrbracket$, $$a_1\cdot\prod_{i=2}^m\min(a_{i-1}+1, a_i)\leq c(\llbracket a_1, a_2, \dots, a_m\rrbracket) \leq a_1\cdot\prod_{i=2}^m\left[a_i+\sum_{j = 1}^{i-1}2a_j\right].$$
\end{theorem}

In addition, we modified a computer algorithm from \cite{sawadali} to count $c(\llbracket a_1, a_2, \dots, a_m\rrbracket)$ for various classes of MV assignments. From this we conjectured that the number of ways to fold MV assignments with equal block sizes is always polynomial. In particular, if we let $S(m,k)$ denote the MV assignment with $m$ blocks, each of size $k$, then we conjecture the following:

\begin{conjecture}
    For any fixed $m$ and sufficiently large $k$, $$c(S(m,k))\sim \frac{e^m}{\sqrt{2\pi m}} k^m.$$
\end{conjecture}

We also tried to look for MV assignments that had a maximum number of foldings. A summary of our data is shown in \cite{stamp1} and led us to the following:

\begin{conjecture}\label{conj:max}
    For $n \geq 8$, the MV assignment that folds in the most ways consists only of blocks of size 1, 2, and 3. Furthermore, there are no two consecutive blocks of size 3.
\end{conjecture}

\begin{figure}
\centering
    \begin{tikzpicture}[scale=0.75]
    \begin{semilogyaxis}[
        enlargelimits=false,
        xlabel={$n$},
        ylabel={$\#$ of ways},
        ymajorgrids=true,
        xmajorgrids=true,
        grid style=dashed,
    ]
    \addplot[
        black,
        only marks,
        mark=*,
        mark size=1.5pt]
    table{scatter_one.dat};
    \addplot[
        black,
        only marks,
        mark=o,
        mark size=1.5pt]
    table{scatter_two.dat};
    \addplot[
        domain=4:38,
        color=black,
        ] { 2^x/(x^(5/4)) };
    \end{semilogyaxis}
\end{tikzpicture}
\caption{Number of ways the maximally foldable assignment folds on the $1\times n$ strip of stamps. The curve $2^n/n^{5/4}$ is overlayed.}\label{fig:graph}
\end{figure}
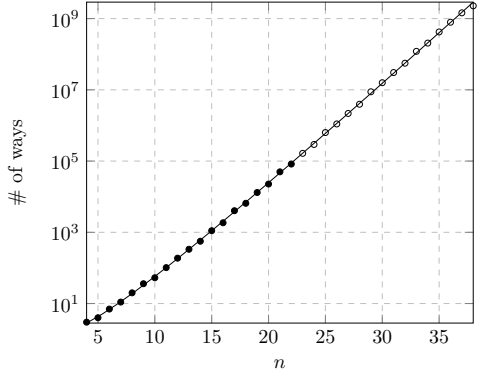

Assuming the above conjecture is correct, we can find the maximally folding assignment and the number of ways it folds for larger $n$. By only considering assignments that satisfy the conditions, we have been able to plot $n$ vs the number of ways the maximally foldable assignment folds for up to $n=33$ and found it to be roughly linear on a logarithmic scale that seems to closely fit $2^n/n^{5/4}$, as shown in Figure~\ref{fig:graph}. 

\begin{conjecture}\label{conj:maxasym}
The maximal number of ways to fold a $1 \times n$ strip approaches $2^n/n^{5/4}$ as $n\to\infty$.
\end{conjecture}

There still remains much to be explored in the fixed MV assignment perspective on the stamp folding problem. 

\begin{openq}
    Is there a closed form for finding the number of ways a specified MV assignment can fold flat? If not, are there any more general closed forms we can generate for interesting patterns?
\end{openq}

\begin{openq}
    Is it possible to compute $c(\mu)$ faster than the algorithm described here?
\end{openq}

An interesting related area of research involves folding the $2\times n$ grid. This problem is somewhat more complicated as the folding cannot be performed recursively, as is the case in the $1\times n$ strip (see \Cref{prop:induction}). In fact, determining if a specific MV assignment on the $2\times n$ has a valid folding is a challenging problem: \cite{Morgan2012MapFB} gives an $O(n^9)$ algorithm. 

\begin{openq}
    Can the insights from the $1\times n$ strip of stamps be used in any way to provide insight towards the $2\times n$ case. For each MV assignment on the $1\times n$ strip, is there a corresponding assignment that folds similarly on the $2\times n$ strip? 
\end{openq}

\acknowledgements
\label{sec:ack}
The authors thank the other participants of the MathILy-EST 2024 REU for numerous conversations and support.

\nocite{*}
\bibliographystyle{abbrvnat}
% use the following instead if you encounter problems 
%\bibliographystyle{alpha}
\bibliography{references}
\label{sec:biblio}

\end{document}